\documentclass[runningheads]{cl2emult}

\usepackage{makeidx}  
\usepackage{graphicx} 
\usepackage{subeqnar} 
\usepackage{multicol} 
\usepackage{cropmark} 
\usepackage{lnp}      
\makeindex            

\begin{document}
\title*{ Diophantine Conditions and Real or Complex Brjuno Functions}
\toctitle{ Diophantine Conditions and Real or Complex Brjuno Functions}
%
%
\titlerunning{ Real or Complex Brjuno Functions}
%
\author{Pierre Moussa\inst{1}
\and Stefano Marmi\inst{2}}
\authorrunning{Pierre Moussa and Stefano Marmi}
%
%
\institute{Service de Physique Th\'eorique, CEA/Saclay,\\
F-91191 Gif sur Yvette cedex, France
\and Dipartimento di Matematica ``U.Dini''
     Universit\`a di Firenze,\\
     Viale Morgagni 67/A, I-50134 Firenze, Italy}

\maketitle              

\begin{abstract}
The continued fraction
expansion of the real number $x=a_0+x_0$,
$a_0\in {\bbbz}$,  is given by $0\le x_n<1$,
$x_{n}^{-1}=a_{n+1}+ x_{n+1}$, $a_{n+1}\in {\bbbn}$,
for $n\ge 0$. The Brjuno function is then
$B(x)=\sum_{n=0}^{\infty}x_0x_1\ldots
x_{n-1}\ln(x_n^{-1})$, and the number $x$  satisfies
the Brjuno diophantine condition whenever  $B(x)$ is bounded.
Invariant circles under a  complex
rotation persist when the map is analytically
perturbed, if  and only if
the rotation number satisfies the Brjuno condition,
and the same holds for invariant circles in the semi-standard and
standard maps cases.
In this lecture, we will review some properties of the Brjuno function,
and give some generalisations  related to  familiar
diophantine conditions. The Brjuno function is highly singular
and takes  value $+\infty$ on a dense set  including
rationals. We present a regularisation leading to
a complex function  holomorphic in the upper half plane.
Its  imaginary part tends to the Brjuno function on
the real axis, the real part remaining bounded, and we also
indicate its transformation under the modular
group.
\end{abstract}

\section{Hamiltonian Chaos and the Standard Map}
The simplest non trivial model for Hamiltonian chaos is a two dimensional
real map, called the ``Standard Map''. It has been introduced more or less
independently by Chirikov and Taylor \cite{chir,esca}.
The occurrence of chaos was discussed by Greene \cite{gree} who
displayed many
numerical results on this model, which describes a simplified version of
the non--linear coupling of two oscillators.
It occurs naturally in many domains of physics, including celestial
mechanics, classical quasiperiodic systems, quantum quasicrystals, adiabatic
response in non--linear mechanics, magnetic toroidal configurations in plasma
physics, non--linear electronic devices, and many others.\par
The Standard Map is a map from the cylinder ${\bbbt}\times {\bbbr}$
to itself, defined as
\begin{equation}
\big(\theta,r\big)\mapsto \big(\theta',r'\big)=
\bigg(\theta+r +{K\over 2\pi}\sin(2\pi\theta)
\ \ ({\rm mod}\ 1)\ ,\ r +{K\over 2\pi}\sin(2\pi\theta)\bigg)\; .
\end{equation}
Note that the second variable can also be taken modulo 1, in which
case we get a map ${\bbbt}^2 \to{\bbbt}^2$. The above map can be
written in two equivalent forms:
\begin{itemize}
\item Hamiltonian form
\begin{equation}
 r'-r={K\over 2\pi}\sin(2\pi\theta)\quad,\quad \theta'-\theta=r'\; ,
\end{equation}
\item Lagrangian form,  where one considers now the twice iterated map,
that is $(\theta,r)\to (\theta',r')\to (\theta'',r'')$, so that
\begin{equation}
\theta''-2\theta'+\theta={K\over 2\pi}\sin(2\pi\theta')\; ,
\end{equation}
and for the $n$-th iterated map one gets
\begin{equation}
\theta_{n+1}-2\theta_n+\theta_{n-1}={K\over 2\pi}\sin(2\pi\theta_n)\; .
\end{equation}
\end{itemize}
This last equation is sometimes called the Frenkel-Kontorova
model \cite{aubr}
which describes equilibrium positions of a chain of material points
placed in a periodic potential and submitted to an harmonic
elastic force between two neighbouring points.

If $K=0$, we get the so-called ``twist-map'', which gives
\begin{equation}
r_{n+1}=r_n=r_0=\rho={\hbox{constant}}\quad,\quad
\theta_{n+1}=\theta_0+n\rho
\ (\hbox{mod}\ 1) \;,
\end{equation}
after $n$ iterations.
This  map is nothing but a rotation of angle
$2\pi\rho$---we say that the rotation
number is $\rho$. The orbits are all transverse to the axis of
the cylinder. They are made of a finite number of points
(and therefore discrete) if
$\rho$ is rational, and they are dense  on transverse circles if
$\rho$ is irrational. The cylinder is sliced along orbits at irrational
values
of $r$, which are intertwined with discrete orbits at rational values of $r$.
The question is : in which sense is such a pattern stable under perturbations,
that is here when $K\neq 0$?

Among the orbits which are dense in a curve
wrapped around the cylinder, particularly interesting are the orbits
which will persist under perturbation, in particular because
they separate the space into
domains which do not communicate. It is known that when $K$ is
large (for example $K>4/3$, see \cite{mather}), such orbits do
not exist, and on the other hand, when $K$ is small,
some of the irrational orbits persist, depending on arithmetical properties
of the rotation number. For perturbed twist maps, we define the rotation number
as $\lim_{n\to\infty} n^{-1}\theta_n$, where
$(\theta_n,r_n)$ is the $n$-th iterated map obtained from (1).

Other kinds of invariant curves  may occur, attached to elliptic periodic
orbits. For instance if $K>0$ is sufficiently small, the point $(1/2,0)$
is an elliptic fixed point. Due to KAM Theorem (see \cite{yoc1} for a review),
there exist homotopically trivial invariant curves on
the cylinder winding around this fixed point, which form the so--called
elliptic islands. We shall not consider here the
problem of such orbits, although there existence is very important in
connection with ergodic theory. Indeed one  expects chaotic behaviour
for $K$ large, but the persistence of elliptic islands could prevent
the map from being ergodic.

\section{The Critical Constants}
For the standard map, we consider now
the homotopically non--trivial invariant curves {\it i.e. }
wrapped around the cylinder. A natural way to look for their existence is to
replace  the angular variable $\theta$ by the  new
variable $\phi\in{\bf T}$
\begin{equation}
\theta=\phi+u(\phi,K,\rho)\quad,\quad
r=\rho+u(\phi,K,\rho)-u(\phi-\rho,K,\rho)\;.
\end{equation}
With the condition $1+{\partial u/\partial\phi}>0$,  it would describe a
curve around the cylinder, on which the map is expressed as
$\phi'=\phi+\rho$,
when $\phi$ describes $\bf T$ for $K$ and $\rho$ fixed. We say that (6)
expresses  on the curve the conjugacy of the map to a rotation. The existence
of a function $u(\phi,K,\rho)$, analytic in the variable $\phi$,
insures the existence of an analytic
invariant curve with rotation number $\rho$. We are interested to determine
the
critical constant $K_c(\rho)$ as being the largest possible value of $K$ for
which such an analytic function $u$ exists. Of course, one could consider
regularity constraints weaker than analyticity, leading to other critical
constants.  We look for a
perturbation expansion of the function $u$, and we follow the notations
of \cite{berge1}.
{}From the standard map we get from (6)
\begin{equation}
u(\phi+\rho, K,\rho)-2u(\phi,K,\rho)+u(\phi-\rho,K,\rho)=
{K\over 2\pi}\sin(2\pi(\phi+u(\phi,K,\rho))\;.
\end{equation}
For $k\ge 1$, we have
\begin{equation}
u^{(k)}(\phi+\rho,\rho)-2 u^{(k)}(\phi,\rho)+
u^{(k)}(\phi-\rho,\rho)= \left.
{1\over 2\pi}\sin\big(2\pi\phi+2\pi u(\phi,K,\rho)\big)\right|_{k-1}\;,
\end{equation}
where in the right hand side, one keeps only the terms of
order $k-1$ in the expansion on powers of $K$. We use now the Fourier series
expansion on $\phi$, that is $u^{(k)}(\phi,\rho)=\sum_{\nu\in{\bf Z}}
u^{(k)}_{\nu}(\rho) e^{2i\pi\nu\phi}$, and we see that the coefficient of
$e^{2i\pi\nu\phi}$ in the left hand side of (8) is $2(\cos(2\pi\nu\rho)-1)
u^{(k)}_{\nu}(\rho)$. Therefore (8) allows a recursive  computation of
the Fourier coefficients $u^{(k)}_{\nu}(\rho)$, and we get for
$u^{(k)}(\phi,\rho)$ expressions as  trigonometric polynomials in $\phi$.
However terms  of the kind
$2(\cos(2\pi\nu\rho)-1)$ occur in the denominators along the steps of the
recursion. Such factors are called ``small divisors'', some of them vanish when
$\rho$ is rational, and may become arbitrarily small when $\nu$ becomes
large, for irrational $\rho$.
Now let $\widetilde{K}_c(\rho)$ be the minimum over $\phi$
of the convergence radius of the expansion
\begin{equation}
u(\phi, K,\rho)=
\sum_{k=1}^{\infty} K^k u^{(k)}(\phi,\rho) \;.
\end{equation}
For $\rho$ rational, (8) cannot be solved, and we set
$\widetilde{K}_c(\rho)=0$. For  irrational values of
$\rho$, Berretti and Gentile \cite{berge2} were able to control
$\widetilde{K}_c(\rho)$ using the Brjuno function
$B(\rho)$ which is a number theoretic function which
will  define in the following
Section 4. More precisely, there exists $C>0$ such that, for any
irrational $\rho$,
\begin{equation}
\left|\ln\left( (\widetilde{K}_c(\rho))^{-1}\right)
-2B(\rho)\right|<C\;.
\end{equation}
The functions $\widetilde{K}_c(\rho)$ and $e^{-2B(\rho)}$ both
vanish on all rationals,
but the previous equation shows that the ratio
$\widetilde{K}_c(\rho)/e^{-2B(\rho)}$ remains uniformly bounded at
every irrationals.
The fact that this ratio remains bounded is in itself amazing, but
it recalls earlier and  now classical results by Yoccoz \cite{yoc2}
on the linearisation of holomorphic maps. We shall see
later that we may have even better results in the framework of
holormorphic maps.

The determination of the radius of convergence in (9) is not the
whole story.  It is  possible that the function $u$ may be
analytically continued for real values of $\phi$ and
for $K>\widetilde{K}_c(\rho)$ real.
Thus we would have another critical constant
$K_c(\rho)$ such that we still have real analytic curves
for $\widetilde{K}_c(\rho)<K<K_c(\rho)$ real.
The numerical results \cite{stark} seem to indicate that this is
indeed the case: see \cite{carlett} for a detailed discussion of
this issue both from the numerical and the analytical points of view,
which also uses results of \cite{russk,gelf}.
The definitive
answer is not known to us today, although we are led to
expect that the function $B(\rho)$ plays
a central role in the determination of $K_c(\rho)$ (see also  Davie
\cite{davi1}).

\section{Complex Analytic Maps} The problems of the critical constant is better
understood in the case of the complex analytic maps. We have already seen
in (10) that the critical constant $\widetilde{K}_c(\rho)$ of the
complexified  version of the standard map is controlled by the
Brjuno function. A simpler example is the ``Semi-standard Map'', which is a
two dimensional complex map on the cylinder, closely related to the standard
map (1) : to get the semi-standard map, just
replace in (1) the sine function
$\sin(2\pi\theta)$ by its positive frequency part $(1/2i)\exp(2i\pi\theta)$.
The procedure to get analytic invariant curves proceeds in a completely similar
way as Equations (6) to (9), and it was proven that in this
case \cite{davi2,marm1}, the critical
constant $K_{\rm ssm}(\rho)$  defined in a same  way as above,
fulfils as in (10)
\begin{equation}
\left|\ln\left( ({K}_{\rm ssm}(\rho))^{-1}\right)
-2B(\rho)\right|<C\;.
\end{equation}
The numerical results (especially the figure 16) in
ref. \cite{marm1} provide
more. Not only the ratio
${K}_{\rm ssm}(\rho))/ e^{-2B(\rho)}$ is bounded on irrationals, but it
is extendable to a continuous function on $[0,1]$, bounded below and
above by positive constants. This result is amazing if one remember that both
${K}_{\rm ssm}(\rho))$ and $ e^{-2B(\rho)}$  vanish at all rationals.
Therefore the Brjuno function $B(\rho)$ is a good model to represent the
singular behavior of
$\ln\big(({K}_{\rm ssm}(\rho))^{-1}\big)$.

The Brjuno function was introduced by Yoccoz \cite{yoc2} in the apparently
simpler problem of the linearisation of complex holomorphic maps around their
fixed points. This is a more than
one century old problem (see \cite{perez} for a nice
review), which we can state as follows. Let $f(z)$ be a holomorphic map such
that $f(0)=0$, $f'(0)=e^{2i\pi\rho}$. Is it possible to conjugate the map
$f$ to its linear part?  This means that we look for a function $h$,
holomorphic in a disk of radius $R_f$, such that
$h(0)=0$, $h'(0)=1$, and $f(h(z))=h(ze^{2i\pi\rho})$.
Note that such a function $h$, if it exists, is unique. In this case,
the function $f$ is said to admit a Siegel disk of radius $R_f$.
The Siegel disk is  a topological disk  with conformal radius $R_f$,
since it is the image through
the  normalised conformal map $h$ of the disk $|z|<R_f$.

We quote now the classical
results  on this question \cite{perez}.  {\it i)} If $\rho$ is rational,
there is no disk, that is
$R_f=0$. {\it ii)} if $\rho$ is irrational and satisfies a (strong)
Liouville condition,
we still have $R_f=0$. {\it iii)} if $\rho$ is a diophantine irrational (see
Section 5 below), then there exists a Siegel disk and $R_f>0$, more precisely
this happens when $B(\rho)$ is finite. {\it iv)} If $B(\rho)=+\infty$, then
that there exist functions $f$ such that $R_f=0$.
Indeed Yoccoz \cite{yoc2} proved the following : define $R(\rho)$ as the
smallest radius of the Siegel disks $R_f$ obtained when $f$ varies in the
compact
family of all univalent maps on the unit disk such that $f(0)=0$, and $f'(0)=
e^{2i\pi\rho}$. Then we have
\begin{equation}
\left|\ln\left((R(\rho))^{-1}\right)-B(\rho)\right|<C\;.
\end{equation}
Now consider the family of quadratic polynomial $P_2(z)=
e^{2i\pi\rho}(z-z^2)$, and call $R_2(\rho)$ the radius of the Siegel disk
associated to it. Observe first  that, through the rescaling
$z \to e^{-2i\pi\rho}R\times z$, then
$P_2$ is transformed in $e^{2i\pi\rho}z-Rz^2$. In the rescaled
variable, we see that $R_2$ is the maximum value of the constant $R$
for which a circle with
conformal radius one is persistant. Therefore $R_2(\rho)$ is the
critical constant adapted to the present case, and this leads to
the analogies between (10), or its equivalent in the  real case,
and (12).

Here  again,  the numerical results (now  the figure 6) in
ref. \cite{marm1} bring some continuity properties.
Not only the ratio
$R_2(\rho))/ e^{-B(\rho)}$ is bounded on irrationals, but it
is extendable to a continuous function on $[0,1]$, bounded below and
above by positive constants.
Therefore the Brjuno function $B(\rho)$ is again a good model to represent the
singular behavior of
$\ln\big((R_2(\rho))^{-1}\big)$.
In our work \cite{marm2}  which started from these observations,
we give arguments which strongly support the conjecture that the ratio
$R_2(\rho))/ e^{-B(\rho)}$
is not only continuous but satifies a H\"older  continuity
condition with exponent $1/2$.
More precisely, the Brjuno function displays the universal  singular
behaviour (up to some H\"older-${1\over2}$ continuous function)
of the critical functions occuring in
small divisors  holomorphic problems in dimension one.

The use of the Brjuno function was somewhat implicit in the work of Buric et
al. \cite{buric}, where they attempted to find representations of the critical
constants by what they called modular smoothing. Singular functions of the same
type occured in MacKay \cite{mackay} in relation to the Brjuno condition.\par
It is nevertheless
useful to recall here briefly one of the steps, called renormalisation,
which plays a special role in Yoccoz's argument,  and will appear
to be crucial in understanding
the fine regularity properties of ratios of the type
$R_2(\rho))/ e^{-B(\rho)}$. For this purpose,  we follow \cite{perez},
and we we consider first a rotation of angle  $2\pi\rho$, with
$0\le \rho<1$, that is $z\to e^{2i\pi\rho}z$, acting
in an open disk of radius $R_{\rho}$ centered at the origin in the complex
plane.  We need an arbitrary point $a$, such that $|a|=R_{\rho}$,
and for simplicity
we take the real point $a=R_{\rho}$. Let $a'=e^{2i\pi\rho}R_{\rho}
$ its image.  Now, consider the angular sector bounded by the lines
$0a$ and $0a'$, namely
$\Delta_{\rho} =\{\;z\;|\; 0\le {\rm arg}\; z<2\pi\rho\;,\;
|z|<R_{\rho}\}$,
the line $0a'$ being excluded.
Consider the orbit  made of the successive iterated points
$z_n,\; n\ge 1$ starting
from $z_0\in \Delta_{\rho}$, and let $z_q$ the first of these points
which also belong to $\Delta_{\rho}$.
The map $z_0\to z_q$ is thus the first return map in the sector.
We have
$(2\pi)^{-1}({\rm arg}\; z_q)=(2\pi)^{-1}({\rm arg}\; z_0)+q\rho -1$.
We now take in the sector the variable
$u$ such that its complex conjugate
$\overline{u}=z^{1\over\rho}$, where
now $u$ belongs to a disk of radius $R_{\rho}^{1\over \rho}$.
For the values $u_0$ and $u_q$, corresponding to $z_0$ and $z_q$,
we have
$(2\pi)^{-1}({\rm arg}\; u_q)=(2\pi)^{-1}({\rm arg}\; u_0)-q +{1\over \rho}
\;({\rm mod}\ 1)\;=(2\pi)^{-1}({\rm arg}\; u_0) +{1\over \rho}$.
The original map which acted in a disk of radius $R_{\rho}$, leads
in the new``renormalised variable'' $u$, to a rotation
with rotation number ${1\over\rho}$, acting in a disk with radius
$R_{(1/\rho)}$, such that
$\ln R_{\rho}=\rho \ln R_{(1/\rho)}$.

This construction extends to the non linear perturbed case, for example
$P_2(\rho,z)= e^{2i\pi\rho}(z-z^2)$,
with a lot of complications. Suppose that there is a Siegel disk
for $P_2(\rho,z)$. In this disk, there are conformal coordinates
on which the maps is exactly a rotation of angle $2\pi\rho$,
and on these coordinates we apply the linear renormalisation.
The problem is then to give an interpretation of the renormalised
coordinates
$u$ which we obtain. It appears that there exist a holomorphic map
in the variable $u$
with rotation number $\rho^{-1}$, which admits a Siegel disk,
with conformal radius  $\widetilde{R}_{(1/\rho)}$ such that
$\ln R_{\rho}-\rho \ln \widetilde{R}_{(1/\rho)}=0$. However, this map
is not a polynomial with degree 2. This led Yoccoz to extend the problem
to the  compact family of univalent map on the unit disk
with rotation number $\rho$, and he  has considered  the minimum
$R(\rho)$ of the radius of the Siegel disk taken
over this family of maps. The result is two modifications to the
relation $\ln R_{\rho}-\rho \ln R_{(1/\rho)}=0$
obtained in the linear case.
First due to the minimisation procedure, the best one could get is
a positive uniform upper bound for this expression instead of zero.
Second, there is a special difficulty when $\rho$ goes to zero. In this case
the Siegel disk is strongly distorted, since there is an other fixed point
which
tends to zero when $\rho$ goes to zero. The comparison between the linear and
the non linear case becomes  unjustified in this limit.
Yoccoz proved that the result is
an additional logarithmic term in the estimate, so that we only get that
$\ln R_{\rho}-\rho \ln R_{(1/\rho)}-\ln \rho$ is bounded. It is therefore
natural to compare the function $-\ln R(\rho)$ (as well as
$-\ln R_2(\rho)$) to the solution of the equation
$B(\rho)-\rho \ln B(1/\rho)+\ln \rho=0$ which we will see,  is nothing else
than the Brjuno function.

\section{Continued Fractions and the Brjuno Function}
We first give a somewhat unusual definition of the continued fraction
expansion sometimes called ``japanese continued fractions'' \cite{schweiger}.
Let $\alpha$ be a fixed real number  such that
${1\over 2}\le \alpha \le 1 $. Then,
given the starting number $x$, the coefficients $a_n$ and $\varepsilon_n$
are  recursively uniquely defined by the conditions
\begin{equation}
x=a_0+\varepsilon_0 x_0,\ \hbox{and}\ \forall n\ge 0,\
x_n^{-1}=a_{n+1}+\varepsilon_{n+1}x_{n+1}\;,
\end{equation}
with $ \forall n\ge 0,\ \alpha -1\le \varepsilon_n x_n <\alpha $.
We define the  modified integer part $[x]_{\alpha}$ and  the modified
fractional
part $\{x\}_{\alpha}$ as follows,
\begin{equation}
[x]_{\alpha}=[x-\alpha+1]_1 \quad\hbox{and}\quad  \{x\}_{\alpha}=
\{x-\alpha+1\}_1 +\alpha -1\;,
\end{equation}
where $[x]_1$ and  $\{x\}_1$ are the  usual integer and
fractional parts of $x$ (so that $0\le \{x\}_1<1$). With these notations,
we can rewrite (13) as
\begin{equation}
a_0=[x]_{\alpha},\ \varepsilon_0 x_0=\{x\}_{\alpha}, \ \hbox{and},
a_{n+1} =[x_n^{-1}]_{\alpha},\ \varepsilon_{n+1} x_{n+1}=
\{x^{-1}_n\}_{\alpha}, \forall n>0 \;.
\end{equation}
Therefore the $x_n$ are generated by iterating the function
$A_{\alpha}(x)=\left|\{x^{-1}\}_{\alpha}\right|$,
that is $ \forall n\ge 0,\ \  x_{n+1}=A_{\alpha}(x_n)=
\left|\{x^{-1}_n\}_{\alpha}\right|=
\left| x_n^{-1}-[x_n^{-1}]_{\alpha} \right|$.
A more detailed description states that the map $A_{\alpha}$ is made of the
following branches
\begin{subeqnarray}
\hbox{branch}\ k^+:\quad A_{\alpha}(x)&=&{1\over x}-k
\quad \hbox{for}\quad{1\over k+\alpha}<x\le{1\over k}
\ts\ , \label{br1a}\\
\hbox{branch}\ k^-:\quad A_{\alpha}(x)&=&k- {1\over x} \quad
\hbox{for}\quad{1\over
k}<x\le{1\over k+\alpha -1}\ts\ . \label{br1b}
\end{subeqnarray}
When ${1\over 2}<\alpha\le 1$, the function $A_{\alpha}$
maps the interval
$[0,\alpha)$ to itself, whereas when $ \alpha={1\over 2}$,
it maps the interval
$[0,\alpha]$ to itself. In both cases, it is convenient to set
$A_{\alpha}(0)=0$,
and we get a map which is infinitely differentiable by pieces, and
the points where it is  not differentiable accumulate to $0$.
Now $x$ and the reduced fraction $p_n/q_n$
admit the following representation
\begin{equation}
x=a_0 + \displaystyle{\varepsilon_0 \over a_1
+ \displaystyle{\varepsilon_1  \over a_2 + \ddots +
\displaystyle{\varepsilon_{n-1} \over a_n + \varepsilon_n x_n}}}
\quad ,\quad
{p_n \over q_n}= a_0 + \displaystyle{\varepsilon_0 \over a_1
+ \displaystyle{\varepsilon_1  \over a_2 + \ddots +
\displaystyle{\varepsilon_{n-1} \over a_n }}}
\ \ .\
\end{equation}
As long as the $x_n$'s do not vanish, we have
\begin{equation}
x={p_n+p_{n-1}\varepsilon_nx_n\over q_n+q_{n-1}\varepsilon_nx_n}
\quad,\quad
x_n=(-\varepsilon_n){p_n-xq_n\over p_{n-1}-xq_{n-1}} \ ,
\end{equation}
and the recursion relations
\begin{subeqnarray}
p_n&=&a_np_{n-1}+\varepsilon_{n-1}q_{n-2}\ ,\ p_0=a_0
\ ,\ p_{-1}=1\ts\ ,\label{reca}\\
q_n&=&a_nq_{n-1}+\varepsilon_{n-1}q_{n-2}\ ,\ q_0=1
\ ,\ q_{-1}=0\ts\ ,\label {recb}
\end{subeqnarray}
so that we get $0<q_0\le q_1<q_2<\ldots<q_n<q_{n+1}<\ldots$.
We also define
\begin{equation}
\beta_n=x_0x_1\cdots x_n=(-1)^n \varepsilon_0 \varepsilon_1\cdots
\varepsilon_n (q_nx-p_n)\; ,
\end{equation}
and we have
\begin{equation}
{1\over 1+\alpha}\le \beta_n q_{n+1}={q_{n+1}\over
q_{n+1}+\varepsilon_{n+1}q_nx_{n+1}}\le{1\over \alpha}\ .
\end{equation}
Now there exist  $\lambda(\alpha)$, with $0<\lambda(\alpha)<1$,
and positive constants $C_1$ and $C_2$, such that\cite{marm2}
\begin{equation}
\beta_n<C_1 \lambda(\alpha)^n\quad,\quad q_n>C_2 \lambda(\alpha)^{-n}\ .
\end{equation}
Indeed we have
\begin{subeqnarray}
\hbox{for} \qquad {\sqrt{5}-1\over 2}<\alpha \le 1\ ,
\quad \lambda(\alpha)&=&
\lambda(1)={\sqrt{5}-1\over 2}=0.618...\ts\ , \label{lamd1a}\\
\hbox{and for}\ {1\over 2}\le\alpha \le{\sqrt{5}-1\over 2}\ ,
\ \lambda(\alpha)&=&\lambda\left({1\over 2}\right)
=\sqrt{2}-1=0.414... \ts\ .  \label{lambdb}
\end{subeqnarray}
When $x_n$=0 for some $n$, and $x_m\neq0$ for $m<n$,
then we have $x=p_n/q_n$  which is rational, and we
say that the fraction stops at order $n$
(with our conventions,  we have $x_m=0, \forall n\ge m$). Conversely,
if $x$ is rational, the continued fraction expansion stops at some
finite order $n$.
For $\alpha=1$, we get the  classical Gauss continued fraction expansion
for which all signs $\varepsilon_n=+1$, and for $\alpha=1/2$, we have the
continued fraction to the nearest integer.
Note that when  $\alpha\neq 1$, the results of equations (23a) and (23b)
are not obvious.  For details, and in particular for the extension to
others values of $\alpha$, with  $0\le\alpha<{1\over2}$,
see \cite{marm2,mouss}.

Given a positive real function $f$ on $(0,1)$, the Brjuno series
$B^{(\alpha)}_f(x)$ is the sum (which can be infinite)
of the series with positive terms
\begin{equation}
B_f^{(\alpha)}(x) =\sum_{n=0}^{\infty}\beta_{n-1}f(x_n)
=f(x_0)+x_0f(x_1)+\ldots+x_0x_1\cdots x_{n-1}f(x_n)+\ldots\;,
\end{equation}
where ${1\over2}\le\alpha\le 1$, and for $k\ge 0$,
$x_k$ defined in Eq. (13) or (15). As mentioned above,
when $x$ is rational, we have $x_n=0$ for some $n$, and we use
$x_m=0$ for $m\ge n$,
The following results are easily obtained from the definitions
\begin{subeqnarray}
B_f^{(\alpha)}(x)&=&B_f^{(\alpha)}(x+1)\ts ,\label{bbb1a}\\
B_f^{(\alpha)}(x)&=&xB_f^{(\alpha)}\left({1\over x}\right)+f(x)
\ \hbox{\ for }0< x< \alpha \ts ,\label{bbb1b}\\
B_f^{(\alpha)}(x)&=&B_f^{(\alpha)}(-x)
\ \hbox{\ for }0< x\le 1-\alpha \ts .\label{bbb1c}
\end{subeqnarray}
In particular, $B_f^{({1\over2})}(x)$ is an even function.
More surprising is the following result \cite{marm2} : in the $\alpha=1$ case,
for $B_f^{(\pm)}(x)={1\over2}(B_f^{(1)}(x)\pm B_f^{(1)}(-x))$
which are the even and odd parts of $B_f^{(1)}(x)$, we have
for $0<x\le {1\over 2}$,
\begin{subeqnarray}
\ B_f^{(-)}(x) &=& {1\over2}\left(f(x)-f(1-x)
-(1-x)f\left({x\over 1-x}\right)\right)
\ts ,\label{bbe1a}\\
B_f^{(+)}(x)&=& xB_f^{(+)}\left({1\over x}\right)
+{1\over 2}G(x) \ts ,\quad \hbox{with} \label{bbe1b}\\
\,G(x) &=&f(x)+f(1-x)+(1-x)f\left({x\over 1-x}\right)
+2xB_f^{(-)}\left({1\over x}\right)\ts .\label{bbe1c}
\end{subeqnarray}
In order to prove the previous equations, we use
Equations (23a--c) and the succession of transformations
$$-x\to 1-x\to {1\over 1-x}\to {1\over 1-x}-1={x\over 1-x}
\to {1-x\over x}={1\over x}-1\to {1\over x}\to x \;,$$
which provides the requested relations between
$B(x)$ and $B(-x)$.%

Now, it is convenient to introduce the following specific
notations:\par
{\it i)} In $B_f^{(\alpha)}$(x),  when $\alpha=1$, we
omit the superscript $(\alpha)$, and when $\alpha={1\over2}$,á
we replace  the superscript $(\alpha)$ by $e$, so that
$B_f(x)=B_f^{(1)}(x)$ and $B_f^e(x)=B_f^{(1/2)}(x)$ respectively.

{\it ii)} We omit the subscript $f$ when
$f(x)=-\ln(x)=\ln(x^{-1})$, so that
$B^{(\alpha)}(x)=B^{(\alpha)}_{-\ln}(x)$, $B(x)=B_{-\ln}^{(1)}(x)$
and $B^e(x)=B_{-\ln}^{(1/2)}(x)$ respectively.
We will call $B(x)$ {\it the Brjuno function},
which has been mentioned above in Sects. 2 and 3.
We have
\begin{equation}
B(x)\! =\!\!\sum_{n=0}^{\infty}\beta_{n-1}\!\ln\!\big(x_n^{-1}\big)
\!\!=\!-\ln(x_0)-x_0\ln(x_1)+\ldots-x_0x_1\cdots x_{n-1}
\ln(x_n)-\ldots\;,
\end{equation}
where the $x_n$ are obtained from (13) using $\alpha=1$ (Gaussian case),
whereas $B^e(x)$ is given by the same equation (27) with
$x_n$ obtained from (13) using $\alpha={1\over2}$ (continued fraction to the
nearest neighbour). Both functions  $B(x)$ and  $B^e(x)$ are $1$-periodic,
and  take value $+\infty$ for $x$ rational.
{}From (25a), the odd part of $B(x)$ is given
for $0\le x\le {1\over 2}$  by
$B_{-}(x)={1\over 2}x(\ln(x^{-1}-1))$,
which is continuous (and even H\"older continuous for any exponent
$\sigma<1$).
Moreover, $B^e(x)$ is even, and it has been proven
\cite{marm2} that the difference
$B^e(x)-B^{+}(x)$ is not only bounded, but continuous, and even
H\"older continuous for exponent ${1\over2}$. This refines a more general
statement \cite{marm2,mouss} which says that
the differences $B^{(\alpha)}_{\rm ln}-B(x)$ are bounded over the irrationals.

The numerical computation of $B(x)$ and $B^e(x)$ is delicate,
due to the instabilities of the continued fraction expansion.
However, it is very easy to compute their values when the
continued fraction expansion is periodic, that is when $x$ is
an irrational quadratic number. This applies to noble numbers,
in which case the $x_n$ are constant after a certain order.

\section{The Brjuno Series and Diophantine conditions}
A real number is said to be a {\it Brjuno number}
if and only if $B(x)$ is finite, and we also say that $x$
satisfies {\it the Brjuno  diophantine  condition}.
Brjuno numbers are irrationals and  real numbers
satisfying the classical diophantine conditions (which we recall below)
are Brjuno numbers.
In \cite{marm2}, we show that for ${1\over 2}\le\alpha\le 1$,
$B^{(\alpha)}(x)$
is finite if and only if $x$ is a Brjuno number. More precisely,
the proof says that for any $\alpha\in [0,{1\over2}]$,
the difference $|B^{(\alpha)}(x)-B(x)|$ is bounded
over irrational
values of $x$. We also show that for $\alpha=1$, the difference
$|B(x)-\sum_{n=0}^{\infty} q_{n}^{-1}\ln(q_{n+1})|$ is bounded over
irrational values of $x$, so that we recover the original definition
of the Brjuno numbers \cite{brjuno}:
$x$ is a Brjuno number if and only if
$\sum_{n=0}^{\infty} q_{n}^{-1}\ln(q_{n+1})$ is bounded
over the irrational. One can see \cite{marm2,mouss} that
such  a definition of the Brjuno numbers does not depend of the
particular value of $\alpha$ used to compute the $q_n$.

We now report the usual definition \cite{yoc1} of the
diophantine conditions : we say that $x$ is an irrational diophantine
number of order $\tau\ge 0$
(and we write $x\in {\rm C}(\tau)$), if there exists $c>0$
such that for any integers $p$ and $q$, such that $q>0$, we have
$\left|x-p/q\right|\ge c  q^{-2-\tau}$. Some classical facts
need to be recalled here \cite{hardy}. First, for any $p$ and $q$
such that $0<q<q_{n+1}$, we have $|qx-p|\ge|q_nx-p_n|$, where
$p_n/q_n$ is the Gaussian reduced fraction to $x$. Therefore, in
order to have
$x\in {\rm C}(\tau)$, it is sufficient to check that for any $n>0$,
$\left|x-p_n/q_n\right|\ge c  q_n^{-2-\tau}$. Second, Liouville's
classical theorem asserts that algebraic numbers of degree $n$
belong to $x\in {\rm C}(n-2)$.
Moreover Roth's theorem shows that all algebraic numbers belong to
${\rm C}(\tau)$, for all $\tau>0$.
Finally, for an arbitrary irrational,
and any $n>0$, we have
$(q_n q_{n+1})^{-1}\le \left|x-{p_n/q_n}\right|\le q_n^{-2}$.
Using (21) for $\alpha=1$, we get an equivalent caracterisation
of the diophantine conditions:
$x\in {\rm C}(\tau)$ if and only if there exists a constant $c>0$
such that
$\beta_n\ge c\,\beta_{n-1}^{1+\tau}$ for any $n>0$.

Now we introduce
for $\nu>0$, the Brjuno series
$B_{\{\nu\}}(x)\equiv B_{x^{-\nu}}(x)$
for the fonctions $f(x)=x^{-\nu}$ (still using $\alpha=1$),
\begin{equation}
B_{\{\nu\}}(x)=
\sum_{n=0}^{\infty}\beta_{n-1}(x)x_n^{-\nu}= \sum_{n=0}^{\infty}
\beta_{n-1}\left({\beta_n\over\beta_{n-1}}\right)^{-\nu}\!\!=
\sum_{n=0}^{\infty}\beta_{n-1}^{1+\nu}(x)\beta_n^{-\nu}(x) \ .
\end{equation}
Using (21) one gets
\begin{equation}
2^{-\nu}\sum_{n=0}^{\infty}q_n^{-1-\nu}|q_nx-p_n|^{-\nu}
\le B_{\{\nu\}}(x) \le
\sum_{n=0}^{\infty}q_n^{-1-\nu}|q_nx-p_n|^{-\nu}\ .
\end{equation}
The series
$B_{\{\nu\}}(x)$ converges if and only if the series
$\sum_{n=0}^{\infty}q_n^{-1-\nu}|q_nx-p_n|^{-\nu}$
converges, that is if the series
$\sum_{n=0}^{\infty}q_n^{-1-2\nu}\left|x-(p_n/q_n)\right|^{-\nu}$
also converges. As a consequence, if
$B_{\{\nu\}}(x)<\infty$, then
$q_n^{-1-2\nu}\left|x-(p_n/q_n)\right|^{-\nu}$
is bounded, and $x\in {\rm C(1/\nu)}$.
Conversely, assume $\tau\ge 0$, and
$x\in {\rm C}(\tau)$,  then we have
\begin{equation}
B_{\{\nu\}}(x)\le c^{-\nu}\,
\sum_{n=0}^{\infty}q_n^{-1+\tau\nu} \ .
\end{equation}
Using bounds in (22), we get the following statement:
If $x\in {\rm C}(\tau)$, then for any $\nu$ such that $\tau<\nu^{-1}$,
$B_{\{\nu\}}(x)<\infty$.
Therefore, there is a relation between the diophantine conditions
${\rm C}(\tau)$, and the convergence of the Brjuno series for
$f(x)=x^{-\nu}$: the set of irrationals $x$  such that
$B_{\{\nu\}}(x)\equiv B_{x^{-\nu}}(x)$ is bounded,
is contained in ${\rm C}(1/\nu)$, and contains
${\rm C}(-\varepsilon+1/\nu)$, for any $0<\varepsilon \le 1/\nu$.
In some sense, the Brjuno conditions is related
to the limiting case $\nu=0$, and in particular, $x\in {\rm C}(\tau)$
for $\tau>0$ implies $B(x)<\infty$, that is  $x$ is a Brjuno number.
Using a more general function $f$, positive
on $(0,1)$, and monotoneously decreasing in the vicinity of zero,
we  can introduce a wide family of conditions
$B_f(x)<+\infty$.  The diophantine conditions obtained will be mainly
governed by the singular behavior of $f$ around zero. A power law behaviour
would simulate the usual conditions, whereas a logarithmic behaviour
would generate a condition similar to the Brjuno condition.
Other interesting examples would be obtained by taking
functions $f$ of the form $x^{-\nu}|\log(x)|^{\mu}$,
$x^{-\nu}|\log(x)|^{\mu} |\log(|\log x|)|^{\sigma}$, and so one.

\section{The Brjuno Operator}
We will introduce  now some functional analysis in order to solve
Equations (25a--c).
For fixed ${1\over 2}\le\alpha\le 1$, let us consider the operator
$T_{(\alpha)}$,  acting on locally Lebesgue integrable
functions $f$  on the real line, which verify
\begin{subeqnarray}
f(x)&=&f(x+1) \;  \hbox{for almost every}\, x\in \bbbr
\ts ,\label{per1a}\\
f(x)&=&f(-x)\;\hbox{ for almost every}\, x\in (0,1-\alpha)
\ts .\label{per1b}
\end{subeqnarray}
The operator is defined by
\begin{equation}
(T_{(\alpha)}f)(x)=x f\left({1\over x}\right)\; ,
\ \ \hbox{if}\ \  x\in (0,\alpha )\ .
\end{equation}
It is understood that the function $T_{(\alpha)}f$ is completed
outside $(0,\alpha)$ by imposing on $T_{(\alpha)} f$ the same parity
and periodicity conditions which
are expressed for $f$ in the above equations (31a--b).
The functional equations (25a--c) can then be written
in the form
\begin{equation}
\left(1-T_{(\alpha)}\right)B_f^{(\alpha)}=f\; .
\end{equation}
This suggest to study the operator
$T_{(\alpha)}$ on the Banach spaces
\begin{equation}
X_{\alpha ,p}= \left\{ f : \bbbr \to \bbbr \mid f \,\hbox{verifies (31a--b)}
\; , \;\; f\in L^p(0,\alpha)\right\}
\end{equation}
endowed with the norm of $L^p(0,\alpha)$, namely
\begin{equation}
||f ||_{\alpha ,p}=\left(\int_0^\alpha |f(x)|^p\, dx \right)^{1/p}\; ,
\end{equation}
for $p\in [1,\infty ]$. Note that one could also  use
$L^p(0,1)$, instead of $L^p(0,\alpha)$, and that if $p<p'$
one has the obvious inclusion $X_{\alpha ,p'}\subset X_{\alpha ,p}$.
If $(1-T_{(\alpha)})$ is invertible in the considered space, then
(25a--c) have a unique solution for $B^{(\alpha)}_f$, provided that the
$f$ in the right hand side of (25b) also belongs to the space.
The invertibility  property  is given by the following theorem, which
states in particular that the spectral radius of $T_{(\alpha)}$ is strictly
smaller than 1.\par\noindent
{\bf Theorem.} $T_{(\alpha)}$ is a linear bounded operator from
$X_{\alpha ,p}$ into itself for all $\alpha\in [{1\over
2},1]$ and for all $p\in [1,\infty ]$. Its spectral radius on
$X_{\alpha ,p}$ is bounded by the constant $\lambda(\alpha)$ of
Equation (22), and therefore $1-T_{(\alpha)}$ is invertible.\par\noindent
For the proof, see
\cite{marm2}. We will just observe here
that the result is immediate in the $p=\infty$ case. Indeed,
\begin{equation}
(T_{(\alpha)}^nf)(x) = \beta_{n-1}(x) f(x_n)=
\beta_{n-1}(x)(f\circ A^n_{\alpha})(x)\; ,
\end{equation}
where the map $A_{\alpha}$ is defined above (see (16a--b)).
Therefore
\begin{equation}
||T_{(\alpha)}^nf||_{\alpha ,\infty}\le
{\rm sup}_x(\beta_{n-1}(x))
||f||_{\alpha ,p}\le
c\lambda(\alpha)^{n-1} ||f||_{\alpha ,p}
\end{equation}
and one gets the theorem (for $p$ infinite) by taking the
$1/n$--th root of both sides.  For the other values of $p$, it
is convenient to make use of  the measure which is
invariant under transformation by the map $A_{\alpha}$, instead
of the Lebesgue measure.
An immediate consequence of the theorem, is that if we take $f(x)=\ln(x)$,
for $0<x<\alpha$, then $f\in X_{\alpha,p}$,  for all {\it finite}
$p$ and therefore we also have $B^{(\alpha)}\in X_{\alpha,p}$ for
all finite $p$.\par
However, we have a stronger property in the $\alpha={1\over 2}$
case. Here, we set again $T_e\equiv T_{(1/2)}$.
In this case, the logarithm belongs also to the set
$X_*\subset X_{\alpha=1/2, p=1}$,
made of even, periodic (with period 1) functions belonging to the so-called
``BMO-space''. In this space,  $f$ has bounded mean oscillation,
more precisely  the following semi-norm $||f||_*$
is bounded, with
\begin{equation}
||f||_*=\sup_I{1\over |I|}\int_I|f-f_I|dx\ ,
\end{equation}
where the mean value of $f$ over $I$ is $f_I= |I|^{-1}\int_I|f-f_I|dx$,
and the sup is taken over all possible intervals $I$ with
length $|I|$ smaller than one.
The BMO space has remarkable properties. First it is contained in all
$L_p$ spaces for $p$ finite, and it contains  the $L_{\infty}$ space,
second it is the  space adequate to describe functions having singular
behaviour
not worse than logarithmic, but around every point in a
dense set of the real line,
and third, it has remarkable properties connected to the harmonic conjugacy
transformation \cite{garn,marm2}.

In \cite{marm2,marm3}, we have shown that $1-T_e\equiv 1-T_{(1/2)}$
is invertible in
$X_*$, and therefore  that $B^e\equiv B^{(1/2)}\in X_*$.
Since $|B^{\alpha}-B^e|$ is bounded for ${1\over2}\le\alpha\le 1$,
we have the unexpected consequence  that all $B^{(\alpha)}$
for ${1\over2}\le\alpha\le 1$, also have  the Bounded Mean
Oscillation property, although it cannot be shown directly
through the properties of $1-T_{(\alpha)}$.

The BMO property obtained for the
Brjuno function in the real case, was one of our motivations to consider the
complexification procedure which we will describe in the last Section of this
paper \cite{marm4}.

\section{Application to H\"older--continuous Functions}
In this Section, we will consider only the case $\alpha={1\over 2}$.
In this case the map $A_{\alpha}\equiv A_{1/2}$ is continuous on the interval
$(0,{1\over2}]$.
The functional equation for the Brjuno function $B^e$ for $\alpha=1/2$ is
\begin{equation}
[(1-T_e)B^e](x)=-\log x\; ,
\end{equation}
for all $x\in (0,1/2)$, complemented with the condition that $B^e$
is even
and periodic.  We suppose that the right hand side of
this equation is pertubed, by an additional term $f$, which is
less singular than the
logarithmic function, and we want to study the singular properties of the
perturbed solution. Since  the equation is linear, we only need to consider
the action on $f$ of $T_e$ and $(1-T_e)^{-1}$, which we will conveniently
recall
the
Brjuno operator ${\bf B_e}$. We will consider even and periodic functions $f$
which are {\it continuous}. It is sufficient to know the value of $f$ on
$[0,1/2]$, so we assume $f\in C^0_{[0,1/2]}$. One can check  that
$Tf$ is also continuous provided we set $Tf(0)=0$. We need now the usual
H\"older's type semi-norms for continuous functions :
let $f\in C^0_{[0,1/2]}$, then we define the
H\"older's $\gamma$-norm as
\begin{equation}
|f|_{\gamma}=
\sup_{0\le x<y\le1/2} {\displaystyle|f(x)-f(y)|\over\displaystyle
|x-y|^{\gamma}}\ ,
\end{equation}
with $0<\gamma\le 1$.
This is a seminorm since it vanishes on constant functions, so that we
introduce the norm:
\begin{equation}
||f||_{\gamma}=A|f|_{\gamma}+B|f|_{\infty}\ ,
\end{equation}
where $|f|_{\infty}$ is the $L_{\infty}$ norm of $f$, and
$A$ and $B$ are positive constants which we can choose arbitrarily,
provided that they do not vanish.
We say that $f\in C^{\gamma}$,
if $f\in C^0_{[0,1/2]}$ and $|f|_{\gamma}$ is finite.
We now have:
\par\noindent {\bf Proposition.} $T_e$ is a bounded operator in $C^{\gamma}$,
for the norm $||f||_{\gamma}$, when $0<\gamma\le 1/2$, provided  $B/A$ is
large enough: if $B/A>(2^{\gamma}-2^{-\gamma})^{-1}$, the norm of $T_e$
corresponding to the norm (41) satisfies
$||T_e||_{\gamma}\le 2^{(2\gamma-1)}\le1$. Therefore
for $0<\gamma<1/2$, $T_e$ is a contraction, and $1-T_e$ is invertible.\par
We need the following Lemma
\par\noindent {\bf Lemma.}\ Let $0< y<x\le 1/2$, and define $x_1$ and $y_1$ by
the
following conditions
\begin{equation}
y={1\over n+y_1}\quad,\quad x={1\over m+x_1}\ ,
\end{equation}
with $n\ge2$ and $ m\ge 2$, and $-1/2\le x_1<   1/2 $ and
$-1/2\le y_1<   1/2$, then we have
\begin{equation}
\left\vert\vert x_1\vert- \vert y_1\vert\right\vert\le
{|x-y|\over|x||y|} \ .
\end{equation}
\noindent{\it Proof of the Lemma.} Since $y<x$, we have
$n-m>x_1-y_1>   -1$. therefore $n\ge m$. Let $n-m=p\ge0$.
We have $x-y=xy(p+y_1-x_1)$. So that we need to prove
$||x_1|-|y_1||\le |p+y_1-x_1|$. This is obvious when $p=0$.
We always have $||x_1|-|y_1||\le 1/2$, so the required inequality also holds
when $p\ge 2$. In the remaining case $p=1$, we set
$\eta= \hbox{sign}(y_1)$ and
$\epsilon=\hbox{sign}(x_1)$, and we need to check that
$||x_1|-|y_1||\le |1+\eta|y_1|-\epsilon|x_1||$. Still because the left hand
side is smaller or equal to $1/2$, this last inequality is not obvious
only when $\eta=-1$ and $\epsilon=+1$. It therefore remains to show that
$||x_1|-|y_1||\le |1-|y_1|-|x_1||$. Setting $u=1/2-|x_1|$ and $v=1/2-|y_1|$,
the last inequality is equivalent to $|1-v/u|\le|1+v/u|$, which is readily
checked since $u/v$ is real and non-negative.\par
\noindent{\it Proof of the Proposition.} Let $0<   y<x\le1/2$, and $x_1$ and
$y_1$  as in the preceding lemma. We have
\begin{subeqnarray}
|T_ef(x)-T_ef(y)|&=&|xf(1/x)-yf(1/y)|=|xf(|x_1|)-yf(|y_1|)|
\ts \label{hol1a}\\
&\le&|x-y||f(|x_1|)|+|y||f(|x_1|)-f(|y_1|)|
\ts\label{hol1b}\\
&\le&|x-y||f|_{\infty}+|y||f|_{\gamma}||x_1|-|y_1||^{\gamma}
\ts \label{hol1c}\\
&\le&|x-y||f|_{\infty}+|f|_{\gamma}{|x-y|^{\gamma}\over
|x|^{\gamma}|y|^{\gamma-1}}
\ts , \label{hol1d}
\end{subeqnarray}
where we have used $f\in C^{\gamma}$, and  the Lemma. Therefore
\begin{subeqnarray}
{|T_ef(x)-T_ef(y)|\over|x-y|^{\gamma}}&\le&
|x-y|^{1-\gamma}|f|_{\infty}+\left({|y|\over|x|}\right)^{\gamma}
|y|^{1-2\gamma}|f|_{\gamma}
\ts \label{hol2a}\\
&\le& (1/2)^{1-\gamma}|f|_{\infty}+
(1/2)^{1-2\gamma}|f|_{\gamma}
\ts , \label{hol2b}
\end{subeqnarray}
since $0<   y<x\le 1/2$, and $\gamma\le 1/2$. For $y=0$,
the right hand side can be replaced
by its first term $(1/2)^{1-\gamma}|f|_{\infty}$,  and
the above inequality extends to the case where $y$ vanishes, so that
\begin{equation}
|T_ef|_{\gamma}\le K_{\gamma}(f)=2^{\gamma-1}|f|_{\infty}+
2^{2\gamma-1}|f|_{\gamma}\ .
\end{equation}
For the norm, we get
\begin{subeqnarray}
||T_ef||_{\gamma}&=&A|T_ef|_{\gamma}+B|T_ef|_{\infty}
\ts \label{hol3a}\\
&\le &2^{2\gamma-1}\left[A|f|_{\gamma}+(2^{-2\gamma} B+2^{-\gamma}A)
|f|_{\infty} \right]
\ts ,\label{hol3b}\\
&\le& 2^{(2\gamma-1)}||f||_{\gamma}
\ts ,\label{hol3c}
\end{subeqnarray}
provided $2^{-2\gamma} B+2^{-\gamma}A\le B$, that is
$A/B\le 2^{\gamma}-2^{-\gamma}$
which completes the proof.
The above proposition has two obvious consequences
\begin{itemize}
\item
Since $C^{\gamma}\subset C^{\gamma'}$ whenever
$\gamma'\le\gamma$, we have
\begin{subeqnarray}
&f&\in C^{\gamma}\quad\hbox{and}\quad\gamma\ge 1/2\quad\quad
\Longrightarrow\quad T_ef\in C^{1/2}
\ts \label{hol4a}\\
&f&\in C^{\gamma} \quad\hbox{and}\quad\gamma\ge \gamma_0\quad,\quad
\gamma_0\le1/2\quad \quad \Longrightarrow\quad T_ef\in C^{\gamma_0}
\ts \label{hol4b}
\end{subeqnarray}
\item
When $\gamma<1/2$, $T_e$ is a contraction on $C^{\gamma}$. Therefore
$1-T_e$ is invertible and ${\bf B_e}=\sum_0^{\infty}
(T_e)^n =(1-T_e)^{-1}$ preserves $C^{\gamma}$,
and we have
\begin{subeqnarray}
&f&\in C^{\gamma}\quad\hbox{and}\quad 0<\gamma< 1/2\quad\quad
\Longrightarrow\quad {\bf B_e}f\in C^{\gamma}
\ts \label{hol5a}\\
&f&\in C^{1/2} \quad\quad \Longrightarrow\quad {\bf B_e}f\in C^{\gamma}
\quad,\quad \forall \gamma \ \hbox{such that}\quad 0<\gamma<1/2\ .
\quad\quad\ts \label{hol5b}
\end{subeqnarray}
\end{itemize}
We have reproduced here the proof of \cite{marm3}, because it is
essentially elementary.
In fact we have a slightly better result for ${\bf B}_e$ than for $T_e$,
as shown in the next proposition, which shows that the $C^{1/2}$ property
is effectively reached. Its proof \cite{marm2,marm3} is too difficult to be
reproduced here. We have
\par\noindent {\bf Proposition.} If $f\in C^{\gamma}$\ , and $\gamma>1/2$,
then ${\bf B}_ef\in C^{1/2}$.\par
The Brjuno function $B^e$ which we have studied in the previous section is
nothing else than ${\bf B}_e\ell$, where $\ell$ is equal to minus the
logarithmic fonction restricted to $[0,1/2]$. When made even and periodic,
this function is not continuous. Suppose that we perturb $\ell$ by a function
with enough regularity properties (for example $C^1$ or $C^{1/2+\epsilon}$),
the change in  ${\bf B}_e\ell$ will be continuous and even in
$C^{1/2}$, that is  H\"older-${1\over2}$ continuous.
In this sense, the `most singular  part' of the Brjuno function
is stable or `universal', roughly speaking modulo H\"older-${1\over2}$
continuous contributions. As noticed at the end of Section 5 above,
we can use either $B^e\equiv B^{(1/2)}$ or $B\equiv B^{(1)}$ since
a similar argument starting
from (26a--c) shows that their difference is also H\"older-${1\over2}$
continuous \cite{marm2}.

This provides a frame to understand the properties of
the critical constants $K$ of the Sections 2 and 3 above. We assume
that the singularity comes from the renormalisation equation
$(1-T^e)K=\ell +f$, and not from additional singular behaviour
coming from $f$ in the right-hand side. If it would exist, such an
additional singular behaviour
would require a further physical interpretation. This argument,
which is usual in the renormalisation analysis of singularities,
was  one of the motivation for our previous work \cite{marm2}.
The renormalisation equation allows naturally to conjecture that
the difference between the Brjuno function $B^e$, and the logarithm
of the various critical constants (multiplied by a suitable coefficient),
is continuous and even H\"older-${1\over2}$ continuous.
As an example, we conjecture  that the ratio of $e^{-B}$
and the radius of the Siegel disk of the quadratic polynomial, is an
H\"older-${1\over2}$ continuous function of the rotation number.
These conjectures are in agreement with the numerical results
\cite{marm1,stark}.

\section{ The Complexification of the Brjuno Function}
We consider here the case $\alpha=1$, and we want to associate to the function
$f$ in $X_{1,2}$, a function $\Phi$, holomorphic in the upper half
plane, such that ${\rm Im}\, \Phi\to B_f$ when z goes to the real axis. Since
${\rm Re}\,\Phi$ is associated to the harmonic conjugate of $f$,
we expect to find better boundedness properties when $f$ has the BMO property.
We will here describe our procedure, and report the results \cite{marm4}.

We associate to $f$ a function $F(z)$ holomorphic in $\bbbc\backslash [0,1]$,
and vanishing at infinity, as follows
\begin{equation}
F(z)={i\over \pi}\int_o^1{f(x)\over x-z}\, dx\ .
\end{equation}
For x real, we have ${\rm Im}\, F(x\pm i\varepsilon)=\pm f(x)$ for $
x\in[0,1]$,
and ${\rm Im}\, F(x)=0$ for $ x\not\in[0,1]$.
We will be particularly interested in the case $f(z)=\ln(z)$, in which
case we get $F(z)=-\pi^{-1}{\rm Li}_2(1/z)$,  where ${\rm Li}_2$ is the
classical dilogaritm function \cite{lewin}.  Now we set
\begin{equation}
\Phi(z)=\lim_{N\to\infty}\sum_{-N}^{+N} F(z+n)\ ,
\end{equation}
and we get a function $\Phi$ holomorphic in the upper-half plane $\bbbh_+$,
periodic with a real period equal to one,
and such that for x real, ${\rm Im}\, F(x\pm i\varepsilon)=\pm f(x)$.
In fact the previous equation defines a pair on functions $\Phi_{\pm}$,
respectively holomorphic in the upper or lower half plane $\bbbh_{\pm}$,
so that the natural frame in which our procedure takes place is the frame
of complex hyperfunctions, which we will not consider here \cite{marm4}.

We consider  now the action of $T$, with $(Tf)(x)=xf(1/x)$ if $0\le x<1$,
$f$ and $Tf$ being complemented using periodicity. Using the above
correspondence $f\mapsto Tf$, a correspondence $F\mapsto TF$ is induced
on holomorphic functions in $\bbbc\backslash [0,1]$, vanishing at infinity.
We get
\begin{equation}
(TF)(z)=-z\sum_{m=1}^{\infty}\left( F\left({1\over z}-m\right)-F(-m)\right)
+ \sum_{m=1}^{\infty}F'(-m)\ .
\end{equation}
In fact, $TF$ is essentially $-z\sum_{m=1}^{\infty}F\left(z^{-1}-m\right)$,
up to an affine  additive correction,
which could be determined by the vanishing condition at infinity.

If $f$ is associated to $F$ as above,
the solution $B_f(x)$ (for $\alpha=1$)  of (25a--c) is associated to the
series
\begin{equation}
{\cal B}_f(z)=\sum_{\bbbz}\sum_{m=0}^{\infty} (T^nF(z)) \:,
\end{equation}
where we use the notation
$\sum_{\bbbz}F$ for
$\Sigma_{n=-\infty}^{+\infty}F(z+n)$ understood as the symmetric summation (51)
to insure convergence.

It is now
interesting to display the link between (53) and the modular group
$GL(2,\bbbz)$.
Let $g=\left(
\begin{array}{cc}
a&b\\
c&d
\end{array}
\right)\in GL(2,\bbbz)$,
which mean $a,b,c,d\in \bbbz$, $\varepsilon_g=ad-bc=\pm1$. To
$g$ we associate the following group action on functions holomorphic
on $\bbbc\backslash [0,1]$, that is $g\mapsto L_gF$, with
\begin{equation}
\big(L_gF)(z)=(a-cz)\left\{
F\left({dz-b\over a-cz}\right) -F\left(-{d\over c}\right) \right\}
-{\varepsilon_g\over c} F'\left(-{d\over c}\right)\;.
\end{equation}
Let ${\cal M_+}\subset GL(2,\bbbz)$ be the  multiplicative
monoid generated by the unit matrix, and the set of matrices
$ \left(
\begin{array}{cc}
0&1\\
1&m
\end{array}
\right)$, for $m\ge 1$ integer.
The monoid ${\cal M_+}\subset GL(2,\bbbz)$  can also be defined as
the set of matrices including identity and the matrices
$\left(
\begin{array}{cc}
a&b\\
c&d
\end{array}
\right)\in GL(2,\bbbz)$ such that first,
$d\ge c\ge a\ge 0$, and second, $d\ge b\ge a\ge 0$.
In $GL(2,\bbbz)$  there is a unique  product decomposition, namely
$\forall g \in GL(2,\bbbz)$  there exist a unique set of three matrices
$k$, $m$ and $h$, with $g=kmh$, and $m\in{\cal M}^+$,
$k\in Z$, where Z is the translation subgroup of matrices
$\left(
\begin{array}{cc}
1&n\\
0&1
\end{array}
\right)$, $n\in \bbbz$, and
$h\in H$, where $h$ is the order eight sugroup of
$GL(2,\bbbz)$ made of the matrices
$\left(
\begin{array}{cc}
\varepsilon&0\\
0&\varepsilon'
\end{array}
\right)$, and
$\left(
\begin{array}{cc}
0&\varepsilon \\
\varepsilon' &0
\end{array}
\right)$,
with $\varepsilon=\pm 1$ and $\varepsilon'=\pm 1$.

Now (53) is rewritten as
\begin{equation}
{\cal B}_f(z)=\sum_{k\in Z}\sum_{g\in {\cal M}^+}
\left(L_{(hg)}F\right)(z)\ .
\end{equation}
The double sum over $g$ and $k$ amounts to a sum over a part of the full
modular
group (here one over eight). The contribution over the seven other possible
parts
would be $- {\cal B}_f(z)$, $\pm {\cal B}_f(z^{-1})$, $\pm {\cal B}_f(-z)$, and
$ {\cal B}_f(-z^{-1})$.

We will now summarize the results:\par\noindent
i) The sums in (55) converge in the open upper half plane as long as $f$ is in
$L_1(0,1)$ which insures that $F$ is holomorphic in $\bbbc\backslash [0,1]$,
and vanishes at infinity.\par\noindent
ii) When $f$ is in $L_p(0,1)$, $p$ finite, then
${\cal B}_f$ is in the Hardy $\bbbh_p$ space.\par\noindent
iii) If $f$ is such that $F$ has bounded real part, then the same holds for
${\cal B}_f$.\par\noindent
iv) For $f(x)=\ln(x)$, and $F(z)=-\pi^{-1}{\rm Li}_2(1/z)$,
we get  the complexified Brjuno function,
${\cal B}$, holomorphic in the upper half plane, vanishing at $+i\infty$.
When $z$ goes to a real number in a non-tangential way, we have the following
limits when $\varepsilon>0$ goes to zero$\,$:
the real part ${\rm Re}\,{\cal B}(x+i\varepsilon)$ has a bounded limit
for any real $x$. This limit is continuous at all irrationals and has
a decreasing jump of $\pi/q$ at each rational $p/q$.
When
$x$ is a Brjuno number,
${\rm Im}\, {\cal B}(x+i\varepsilon)$ goes to the Brjuno function $B(x)$.

The limit properties of the complex Brjuno function  on the real axis are
characteristic of functions $f$ having a logarithmic singularity aroud zero.
Although the boundedness of the real part reminds the BMO property of the
real Brjuno function, it  is in fact  a stronger property. This is one
more remarkable feature of this function.
We are convinced that the interpretation of the properties of the
Brjuno function in terms of the modular group is promising.
On the other hand, we  can hope to find an interpretation
of the complexified version of rotation numbers  analogous
with the usual interpretation of complex frequencies in terms of
damped oscillations, but this is another story.

\subsubsection*{Acknowledgements.}
The first author thanks the CNRS and the organizors of the
\'Ecole Thematique held at La Chapelle des Bois, with special thanks
for Michel Planat.
This work begun during a visit of the second
author at the S.Ph.T.--CEA Saclay and at the Dept. of Mathematics of
Orsay during the academic year 1993--1994. This research has been
supported by the  CNR, CNRS, INFN, MURST and a EEC grant.
We thank J.-C. Yoccoz for his constant help and warm encouragements,
and also for allowing us to report here results obtained in
collaboration with him.

\clearpage
\addcontentsline{toc}{section}{Index}
\flushbottom
\printindex


\begin{thebibliography}{7}
%
\addcontentsline{toc}{section}{References}

\bibitem{chir} Chirikov, ~B. (1979) A universal instability of many-
dimensional oscillator systems. Phys. Reports,
{\bf52}, 263--279

\bibitem{esca} Escande, D. (1985) Stochasticity in classical Hamiltonian
systems: universal aspects.
Phys. Reports,á  {\bf121}, 165--261

\bibitem{gree} Greene, J. M. (1979)  A method for determining a stochastic
transition. J. Math. Phys. {\bf 20}, 1183--1201

\bibitem{aubr} Aubry S. and Le Daeron, P. (1983) The discrete
Frenkel--Kontorova model and its extensions. Physica {\bf 8D}, 381--422

\bibitem{mather} Mather J. N. (1984) Non existence of invariant circles.
Ergod. Theor. and Dynam. Sys. {\bf 4}, 301--309


\bibitem{yoc1} Yoccoz J.-C. (1992) An introduction to small divisors
problems, in: From Number Theory to Physics, Waldschmidt M., Moussa P.,
Luck J.-M.,  and Itzykson C. editors,  Springer-Verlag, Berlin,
pp. 659--679

\bibitem{berge1} Berretti A. and Gentile G. (1998) Scaling properties of the
radius of convergence of the Lindstedt series : the standard map. University
of Roma, Italy, preprint

\bibitem{berge2} Berretti A. and Gentile G. (1998) Bryuno function and
the standard map.  University of Roma, Italy, preprint


\bibitem{yoc2} Yoccoz J.-C. (1995) Th\'eor\`eme de Siegel, nombres de
Bruno et polyn\^omes quadratiques. Ast\'erisque, {\bf 231}, 3--88,
(appeared first as a preprint in 1987).


\bibitem{stark} Marmi S. and Stark J. (1992) On the standard map critical
function. Nonlinearity {\bf 5}, 743--761

\bibitem{carlett} Carletti T. and Laskar J. (1999)
Scaling law in the standard map critical function, interpolating
hamiltonian and frequency analysis. (Preprint, Bureau des Longitudes,
Paris, in preparation)

\bibitem{russk} Treshev D. and Zubelevitch O. (1998) Invariant tori in
Hamiltonian systems with two degrees of freedom in a neighborhood of a
resonance. Regular and Chaotic dynamics, {\bf 3}, 73--81

\bibitem{gelf} Gelfreich G. V. (1999) A proof of exponentially small
transversality of the separatrices for the standard map.
Commun. Math. Phys. {\bf 201}, 155--216

\bibitem{davi1} Davie A. M. (1995) Renormalisation for analytic area preserving
maps. University of Edinburgh preprint

\bibitem{davi2} Davie A. M. (1994) the critical function for the semistandard
map.  Nonlinearity {\bf 7}, 219--229

\bibitem{marm1} Marmi S.  (1990) Critical functions for complex analytic maps
function. J. Phys. A: Math.Gen. {\bf 23}, 3447--3474

\bibitem{perez} Perez-Marco R. (1992) Solution compl\`ete du probl\`eme
de Siegel de lin\'earisation d'une application holomorphe autour d'un
point fixe. S\'eminaire Bourbaki nr.753, Ast\'erisque, {\bf 206}, 273--310

\bibitem{marm2} Marmi S.,  Moussa P., and Yoccoz J.-C. (1997)
The Brjuno function and their regularity properties. Commun. Math. Phys.
{\bf 186}, 265-293

\bibitem{buric} Buric N., Percival I. C., and Vivaldi F. (1990)
Critical function and modular smoothing, Nonlinearity {\bf 3}, 21--37

\bibitem{mackay} MacKay R. S. (1988) Exact results for an approximate
renormalisation scheme and some predictions for the breakup of invariant
tori, Physica {\bf 33D}, 240--265, and Erratum (1989) Physica
{\bf 36D}, 358--265

\bibitem{schweiger} Schweiger F. (1995) Ergodic theory of fibered
systems and metric number thory, Clarendon Press, Oxford,

\bibitem{mouss} Moussa P., Cassa A., and Marmi S. (1999)
Continued fractions and Brjuno functions, J. Comput. Appl. Math.
{\bf 105} 403--415

\bibitem{brjuno} Brjuno. A. D. (1971) Analytical form of
differential equations, Trans. Moscow Math. Soc.
{\bf25} 131--288, and, (1972), {\bf 26} 199--239

\bibitem{hardy} Hardy G. H., and Wright E. M. (1938) An introduction
to the theory of numbers, Clarendon Press, Oxford, chapter 11,
fifth edition 1979

\bibitem{garn} Garnett J. B. (1981) Bounded Analytic functions,
Academic Press, New York.

\bibitem{marm3} Marmi S.,  Moussa P., and Yoccoz J.-C. (1995)
D\'eveloppements en fractions continues, fonctions de Brjuno et espaces BMO,
CEA/Saclay, Note CEA-N-2788

\bibitem{marm4} Marmi S.,  Moussa P., and Yoccoz J.-C. (1999)
Complex Brjuno functions, Preprint SPhT/CEA Saclay T99/066, 71 p.


\bibitem{lewin} Lewin L, (1981) Polylogarithms and Associated Functions,
Elsevier North Holland, New York.

\end{thebibliography}
\end{document}